\documentclass[reqno,11pt]{amsart} 
\usepackage{amsmath}
\usepackage{amssymb}
\usepackage{amsthm}

\newcommand\NoBlackBoxes{\global\overfullrule0pt}
\NoBlackBoxes
\theoremstyle{plain}

\begin{document}

\title{Concentration functions and entropy bounds \\
for discrete log-concave distributions
}

\author{Sergey G. Bobkov$^{1}$}
\thanks{1) 
School of Mathematics, University of Minnesota, Minneapolis, MN 55455 USA,
bobkov@math.umn.edu
}

\author{Arnaud Marsiglietti$^{2}$}
\thanks{2)
University of Florida, Department of Mathematics,
Gainesville, FL 32611 USA,
a.marsiglietti@ufl.edu
}

\author{James Melbourne$^{3}$}
\thanks{3)
Department of Electrical and Computer Engineering, 
University of Minnesota, Minneapolis, MN 55455 USA,
jamescmelbourne@gmail.com
}

\subjclass[2010]
{Primary} 

\begin{abstract}
Two-sided bounds are explored for concentration functions and 
R\'enyi entropies in the class of discrete log-concave probability 
distributions. They are used to derive certain variants of the
entropy power inequalities.
\end{abstract}

\maketitle

\def\theequation{\thesection.\arabic{equation}}
\def\E{{\mathbb E}}
\def\R{{\mathbb R}}
\def\C{{\mathbb C}}
\def\P{{\mathbb P}}
\def\Z{{\mathbb Z}}
\def\N{{\mathbb N}}

\def\W{{\rm W}}
\def\T{{\rm T}}
\def\L{{\rm L}}
\def\Re{{\rm Re}}
\def\Im{{\rm Im}}
\def\Tr{{\rm Tr}}

\def\M{{\cal M}}
\def\Var{{\rm Var}}
\def\Ent{{\rm Ent}}
\def\O{{\rm Osc}_\mu}

\def\ep{\varepsilon}
\def\phi{\varphi}
\def\vp{\varphi}

\def\be{\begin{equation}}
\def\en{\end{equation}}
\def\bee{\begin{eqnarray*}}
\def\ene{\end{eqnarray*}}

\vskip5mm
\section{{\bf Introduction}}
\setcounter{equation}{0}

\vskip2mm
\noindent
Given a random variable $X$, its concentration function is defined by
\be
Q(X;\lambda) \, = \, \sup_x \, \P\{x \leq X \leq x+\lambda\}, \qquad \lambda \geq 0.
\en
Of a large interest is also the particular value
\be
M(X) \, = \, Q(X;0) \, = \, \sup_x \, \P\{X = x\}.
\en
These quantities are related to other important characteristics of 
probability distributions such as the Shannon and R\'enyi entropies.
Our original goal was to find effective two-sided bounds on $Q(X,\lambda)$
in the class of discrete log-concave distributions and thus to
explore a number of similarities with the well studied continuous setting.

Let us recall that an integer-valued random variable $X$ (also
called discrete) is said to have a (discrete) log-concave distribution, if 
its probability function $f(k) = \P\{X=k\}$ has an integer 
supporting interval ${\rm supp}(f) = \{k \in \Z: f(k) > 0\}$, and
\be
f(k)^2 \geq f(k-1) f(k+1) \quad {\rm for \ all} \ \ k \in \Z.
\en
Many classical discrete distributions belong to this class: 
discrete uniform, Bernoulli, binomial and convolutions of Bernoulli 
distributions with arbitrary parameters, Poisson, geometric, negative binomial, 
etc. (cf. \cite{J-G} and references therein). It is therefore interesting to know 
how basic characteristics of such distributions are connected with each other:
variance, moments, concentration function, entropies, in analogy 
with usual log-concave distributions in the continuous setting (for which 
the densities with respect to the linear Lebesgue measure are log-concave).
While a number of challenging questions about this class are still open,
in this paper we develop several techniques (reduction to the continuous 
setting, rearrangement) and describe some of the results 
in this direction. In particular, we prove:

\vskip5mm
{\bf Theorem 1.1.} {\sl If the random variable $X$ has a discrete log-concave
distribution, then
\be
\frac{1}{\sqrt{1 + 12\,\Var(X)}} \leq M(X) \leq \frac{2}{\sqrt{1 + 4\,\Var(X)}}.
\en
Moreover, if the distribution of $X$ is symmetric about a point, then
the above upper bound may be sharpened to
\be
M(X) \leq \frac{1}{\sqrt{1 + 2\,\Var(X)}}.
\en
}

\vskip5mm
The inequality (1.5) becomes an equality for two-sided geometric 
distributions, that is, for densities $f(k) = c_p\,p^{|k|}$, $k \in \Z$,
with an arbitrary value $p \in (0,1)$, where $c_p = \frac{1-p}{1+p}$ 
is a normalizing constant. As for the lower bound in (1.4), it does not
need any log-concavity assumption; here the factor
$12$ in front of the variance of $X$ is optimal, as can be seen on 
the example of discrete uniform distributions.

Relations similar to (1.4)-(1.5) with involved parameter $\lambda$ may 
also be stated for the corresponding concentration functions (see Section 8).
Let us however turn to information-the\-ore\-tic applications of Theorem 1.1 
such as entropy power inequalities (EPI's). 

If the random variable $X$ has an absolutely continuous distribution with 
density $f(x)$ with respect to the Lebesgue measure (the continuous setting), 
the differential R\'enyi entropy power of a given order $\alpha > 0$, $\alpha \neq 1$, 
is defined by
\be
N_\alpha(X) = 
\Big(\int_{-\infty}^\infty f(x)^\alpha\,dx\Big)^{-\frac{2}{\alpha - 1}},
\en
while the limit case $N(X) = \lim_{\alpha \rightarrow 1} N_\alpha(X)$
represents the Shannon differential entropy power. As is well-known, 
this functional is subadditive on convolutions, that is, it
satisfies a fundamental EPI
\be
N(S_n) \geq N(X_1) + \dots + N(X_n),
\en
where $S_n = X_1 + \dots + X_n$ is the sum of independent 
continuous random variables. A more general relation such as
\be
N_\alpha(S_n) \, \geq \, c_\alpha \big(N_\alpha(X_1) + \dots + N_\alpha(X_n)\big)
\en
with $c_\alpha = \alpha^{\frac{1}{\alpha - 1}}$, $\alpha>1$, is also true
for R\'enyi entropy powers. We refer an interested reader to \cite{B-G2},
cf. also \cite{M-B}, \cite{R-S}, \cite{B-M}, \cite{L}, \cite{M-M}, and 
\cite{L-M-M} for other variants and extensions of (1.7).

If $X$ takes only integer values with density $f(k)$ with respect to the 
counting measure, the R\'enyi entropy power is defined similarly to (1.6) as
\be
N_\alpha(X) = e^{2H_\alpha(X)} = 
\bigg(\sum_{k \in \Z} f(k)^\alpha\bigg)^{-\frac{2}{\alpha - 1}}.
\en
Here
$$
H_\alpha(X) \, = \,
- \frac{1}{\alpha - 1}\, \log \sum_{k \in \Z} f(k)^\alpha
$$
is the classical R\'enyi entropy, with the limit case
$$
H_1(X) = H(X) = - \sum_{k \in \Z} f(k)\,\log f(k).
$$

Both in the discrete and continuous settings, these entropies are monotone
in the sense that they may only increase when adding an independent summand 
to a given random variable. However, while the EPI (1.7) quantifies 
this property in the continuous setting, not much is known so far about 
the discrete random variables. The inequality (1.7) has been verified to be true 
in \cite{H-V} for the family of the (symmetric) binomial distributions 
(cf. also \cite{J-Y}, \cite{M-W-W}). In the general discrete case, 
it was shown in \cite{H-A-T} that
$$
H(X_1 + X_2) \, \geq \, \frac{1}{2}\,H(X_1) + \frac{1}{2}\,H(X_2) +
g\big(H(X_1),H(X_2)\big)
$$
for some positive doubly-increasing function $g$ on $\R^2_+$ such that
$g(x_1,x_2) \rightarrow \frac{1}{8}$ as $x_1,x_2 \rightarrow \infty$.

Anyhow, as far as we know, the inequality (1.7) is no longer true
in the general discrete case. It seems
the nature of the discrete entropy power is different, since for example
necessarily $N_\alpha(X) \geq 1$. The functionals
$$
\Delta_\alpha(X) = N_\alpha(X) - 1,
$$
reflecting the size of the variance more accurately (in analogy with 
the ``continuous" entropy power), seem to be more appropriate and therefore
more suitable from the point of view of EPI's. 
This can be seen from the following assertion.

\vskip5mm
{\bf Theorem 1.2.} {\sl If the independent random variables $X_k$ $(1 \leq k \leq n)$ 
have symmetric discrete log-concave distributions, then
\be
\frac{1}{c_\alpha}\, \sum_{k=1}^n \Delta_\alpha(X_k) \, \leq \,
\Delta_\alpha(S_n) \, \leq \, c_\alpha \sum_{k=1}^n \Delta_\alpha(X_k)
\en
for any $\alpha>1$ with $c_\alpha = \frac{2(3\alpha - 1)}{\alpha - 1}$.
Moreover, in the case $1 < \alpha \leq 2$, the symmetry assumption may be
dropped, while the constant may be improved to the value
$c_\alpha = \frac{3\alpha - 1}{\alpha - 1}$.
}

\vskip5mm
In this connection, let us mention a recent paper \cite{M-R} dealing
with similar problems. In particular, it is shown there that for 
(not necessarily symmetric) Bernoulli summands $X_k$ and all $\alpha \geq 2$, 
the left inequality in (1.10) remains to hold with
$c_\alpha = \frac{6(\alpha-1)}{\alpha} \leq 6$.
It was also observed that such a lower bound is no longer
true for $\alpha=1$ with any fixed positive constant in place of $c_\alpha$.
It would be interesting to explore whether or not
it is possible to remove the symmetry assumption in Theorem 1.2 for the
whole range of $\alpha$'s.

One may complement these relations with similar ones for the entropy powers
in analogy with (1.8). To this aim, we involve
an additional condition on the variances of the summands.

\vskip5mm
{\bf Theorem 1.3.} {\sl If the random variables $X_k$ have
discrete log-concave distributions, and $\Var(X_k) \geq \sigma^2 > 0$ 
for all $k \leq n$, then for any $\alpha \geq 1$,
\be
N_\alpha(S_n) \, \geq \, \frac{1}{c_\sigma}\, \sum_{k=1}^n N_\alpha(X_k),
\en
where $c_\sigma = 2\pi e\,(1 + \frac{1}{12 \sigma^2})$.
In fact, without any condition on the variances,
this lower bound may be reversed to the form
\be
N_\alpha(S_n) \, \leq \, -\frac{\pi e}{6}\,(3n-1) +
2\pi e \sum_{k=1}^n N_\alpha(X_k).
\en
}

\vskip5mm
The paper is organized as follows. We start with general upper bounds
for the R\'enyi entropy powers in terms of variance and recall known
results in the continuous setting (including an important theorem due to
Moriguti). Such results are used in Section 3 to derive similar upper bounds for
discrete random variables. In Section 4 we turn to the notion of discrete 
log-concavity and first recall its relationship with purely algebraic problems.
Sections 5-6 connects discrete log-concave distributions with
usual (continuous) log-concave measures on the real line. 
In this way, one may derive a number of interesting relations in the 
discrete setting, although not always with sharp constants.
As a sharpening approach, in Section 7 we develop rearrangement arguments,
which allow to complete the proof of Theorem 1.1. A more general form
of this theorem is considered in Section 8 in terms of concentration functions.
Proof of Theorems 1.2-1.3 is postponed to Section 9, and in the last
Section 10 we conclude the exposition with remarks on Bernoulli sums.

\vskip5mm
\section{{\bf Maximum to R\'enyi Entropy Subject to Variance Constraint}}
\setcounter{equation}{0}

\vskip2mm
\noindent
In the discrete case, the $M$-functional (1.2) may be viewed as a member
in the hierarchy of R\'enyi entropies. More precisely, letting 
$\alpha \rightarrow \infty$, the definition (1.9) leads to the identity
\be
N_\infty(X) = M(X)^{-2}.
\en

\vskip5mm
{\bf Convention.} In the continuous case, if a random variable $X$
has density $f(x)$, put
$$
M(X) \, = \, {\rm ess\,sup}_x\, f(x).
$$

\vskip2mm
Thus, we use the same notation for two formally different objects, in analogy 
with $N_\alpha$. As a consequence, the formula (2.1) remains true both 
in the continuous and discrete setting, as follows from the definition (1.6).

First let us state one elementary general relation connecting the $M$-functional
to variance.

\vskip5mm
{\bf Proposition 2.1.} {\sl Given a continuous random variable $X$ with 
a fixed variance, the $M$-functional is minimized for the uniform 
distribution on a finite interval. Equivalently,
\be
M^2(X)\, \Var(X) \geq \frac{1}{12}.
\en
The equality here is attained if and only if $X$ has a uniform distribution
on a finite interval.
}

\vskip5mm
This relation is well-known, and here we recall a simple argument. One may assume 
that $X$ has a finite second moment with $M(X) = 1$ (by homogeneity). Then 
the non-negative function $u(x) = \P\{|X - \E X| \geq x\}$ is Lipschitz 
(therefore absolutely continuous) and satisfies $u(0) = 1$, $u'(x) \geq -2$ a.e.,
so that $u(x) \geq 1 - 2x$ for all $x \geq 0$. This gives
$$
\Var(X) \, = \, 2 \int_0^\infty x u(x)\,dx \, \geq \, 2 \int_0^{1/2} x u(x)\,dx 
 \, \geq \, 2 \int_0^{1/2} x (1-2x)\,dx \, = \, \frac{1}{12}.
$$
\qed

Now, let us restate (2.2) as a homogeneous inequality (with respect to $X$)
\be
N_\infty(X) \leq 12\,\Var(X).
\en
Importantly, this relation may be extended to R\'enyi entropies of all 
orders $\alpha \geq 1$. The problem of maximization of $N_\alpha(X)$ when
the variance is fixed was considered and solved by Moriguti in 1952
(cf. also \cite{K} and \cite{C-H-V} for a multidimensional extension).

\vskip5mm
{\bf Proposition 2.2} (\cite{M}).
{\sl Let $1 < \alpha < \infty$. Given a continuous random variable $X$ with 
a fixed variance, the R\'enyi entropy power $N_\alpha(X)$ is maximized for 
the distribution whose density $f_\alpha$ is supported on the interval $(-1,1)$ 
and is proportional there to $(1 - x^2)^{\frac{1}{\alpha - 1}}$. Equivalently,
\be
N_\alpha(X) \leq A_\alpha \Var(X),
\en
where the constant $A_\alpha$ corresponds to $f_\alpha$.
}

\vskip5mm
The inequality (2.4) is affine invariant, so it is equivalent
to the formally weaker relation $N_\alpha(X) \leq A_\alpha\,\E X^2$
which was actually considered in \cite{M}. 

Thus, the extremal density has the form
$f_\alpha(x) = c_\alpha (1 - x^2)^{\frac{1}{\alpha - 1}}$, $|x| < 1$,
in which the normalizing constant is given by
\be
c_\alpha = \frac{1}{B(\frac{\alpha}{\alpha - 1},\frac{1}{2})} =
\frac{\Gamma\big(\frac{3\alpha-1}{2(\alpha - 1)}\big)}{\Gamma(\frac{\alpha}{\alpha - 1})\,
\Gamma(\frac{1}{2})}.
\en
In this case, as was already noted in \cite{M},
\be
\Var(X) = \frac{\alpha - 1}{3\alpha - 1}.
\en

Let us compute the constant $A_\alpha$. Putting 
$\beta = \frac{2\alpha - 1}{\alpha}$, we have
$$
\int_{-\infty}^\infty f_\alpha(x)^\alpha\,dx \, = \, c_\alpha^\alpha
\int_{-1}^1 (1 - x^2)^{\frac{\alpha}{\alpha - 1}}\,dx \, = \,c_\alpha^\alpha
\int_{-1}^1 (1 - x^2)^{\frac{1}{\beta - 1}}\,dx \, = \,
\frac{c_\alpha^\alpha}{c_\beta}.
$$
Hence, according to the definition (1.4) of the R\'enyi entropy power,
\be
A_\alpha \, = \, \frac{N_\alpha(X)}{\Var(X)} \, = \, \frac{1}{\Var(X)}\,
\Big(\int_{-\infty}^\infty f_\alpha(x)^\alpha\,dx\Big)^{-\frac{2}{\alpha - 1}}
 \, = \, \frac{3\alpha - 1}{\alpha - 1}\
\Big(\frac{c_\beta}{c_\alpha^\alpha}\Big)^{\frac{2}{\alpha - 1}}.
\en
For example,
$
A_2 = \frac{125}{9} \sim 13.888...
$
Although the expression in (2.7) is rather complicated, one can show that
$A_\alpha \rightarrow 2\pi e \sim 17.079$ as $\alpha \rightarrow 1$
(by Stirling's formula), while $A_\alpha \rightarrow 12$ as 
$\alpha \rightarrow \infty$, so that Proposition 2.2 includes (2.3). 
This can also be seen by noting that
the extremal distribution approaches the uniform distribution on the interval 
$(-1,1)$ for large values of $\alpha$, while after a linear transformation 
it approaches the standard normal distribution for $\alpha$ approaching 1.
In the latter case, we arrive at another well-known relation
\be
N(X) \leq 2\pi e\,\Var(X),
\en
where an equality is attained for all non-degenerate normal laws.
Since the function $\alpha \rightarrow N_\alpha(X)$ is non-increasing,
we also conclude that $A_\alpha$ is a decreasing function in $\alpha$.
In particular, (2.8) yields a slightly weaker variant of (2.4) with
a universal constant, namely
\be
N_\alpha(X) \leq 2\pi e\,\Var(X).
\en
It holds for all $\alpha \geq 1$, and an equality is attained for
$\alpha=1$ and all normal laws.

\vskip5mm
\section{{\bf Discrete Case: Bounds on R\'enyi Entropy via Variance}}
\setcounter{equation}{0}

\vskip2mm
\noindent
In order to derive similar relations in the discrete case, one may apply 
(2.4) to the random variable $\widetilde X = X+U$, where $U$ is 
independent of $X$ and has a uniform distribution on the interval 
$(-\frac{1}{2},\frac{1}{2})$. In terms of the probability function 
$f(k) = \P\{X=k\}$, $\widetilde X$ has density
$$
\tilde f(x) \, = \,
\sum_{k \in \Z} f(k)\,1_{(k - \frac{1}{2},\, k + \frac{1}{2})}(x), \quad x \in \R.
$$
It follows that
\bee
h_\alpha(\tilde X) 
 & \equiv &
- \frac{1}{\alpha - 1}\, \log \int_{-\infty}^\infty \tilde f(x)^\alpha\,dx \\
 & = &
-\frac{1}{\alpha - 1}\, \log\, \sum_{k \in \Z} f(k)^\alpha \ = \ H_\alpha(X),
\ene
and therefore $N_\alpha(\widetilde X) = N_\alpha(X)$. Since
$\Var(\widetilde X) = \frac{1}{12} + \Var(X)$, from (2.4) we therefore obtain:

\vskip5mm
{\bf Proposition 3.1.} {\sl Let $1 \leq \alpha \leq \infty$. For any
integer valued random variable $X$ having finite variance,
\be
N_\alpha(X) \leq A_\alpha\, \Big(\frac{1}{12} + \Var(X)\Big),
\en
where the constant $A_\alpha$ is described in $(2.5)$ and $(2.7)$. In particular,
$$
N_\alpha(X) \leq 2 \pi e\, \Big(\frac{1}{12} + \Var(X)\Big).
$$
}

\vskip2mm
For $\alpha = \infty$, we have $A_\infty = 12$, and (3.1) yields:

\vskip5mm
{\bf Corollary 3.2.} {\sl For any
integer valued random variable $X$ having finite variance,
\be
1 \leq N_\infty(X) \leq 1 + 12\,\Var(X).
\en
}

\vskip2mm
In view of (2.1), the above upper bound is exactly the lower bound (1.4) in 
Theorem 1.1. Note that it may also be obtained with a similar argument on the basis 
of Proposition 2.1, thus without referring to Moriguti's theorem
(cf. Proposition 8.2 below).

Since in general $N_\alpha(X) \geq 1$, while the variance may take any prescribed 
value, the constant $\frac{1}{12}$ may not be removed from (3.1). 
Nevertheless, one may ask the following question: Is it possible to replace
$\frac{1}{12}$ in (3.1) with $1/A_\alpha$ at the expense of an additional factor
in front of the variance in analogy with (2.4)? The answer is affirmative in some 
sense for $\alpha>1$, if we allow the factor depend on $\alpha$.

Indeed, generalizing the previous argument, let us apply (2.4) to random 
variables of the form $\tilde X = X + U$, assuming that $U$ is independent 
of $X$ and has density $g(x)$ supported on the unit interval. Then, 
$\tilde X$ has density
$$
\tilde f(x) \, = \, f(k)\, g\Big(x - k + \frac{1}{2}\Big), \qquad
k - \frac{1}{2} < x < k + \frac{1}{2}, \ \ \ 
k \in \Z.
$$
If $\alpha>1$, it follows that
$$
\int_{-\infty}^\infty \tilde f(x)^\alpha\,dx \, = \,
\sum_{k \in \Z} f(k)^\alpha \int_{-\infty}^\infty g(x)^\alpha\,dx,
$$
and therefore
$$
N_\alpha(\tilde X) \, = \, N_\alpha(X)\, N_\alpha(U).
$$
Since $\Var(\tilde X) = \Var(X) + \Var(U)$, 
an application of (3.2) yields
$$
N_\alpha(X)\, N_\alpha(U) \leq A_\alpha\,\big(\Var(X) + \Var(U)\big),
$$
that is,
$$
N_\alpha(X)\, \leq \, \frac{A_\alpha \Var(U)}{N_\alpha(U)}
 + \frac{A_\alpha}{N_\alpha(U)} \Var(X).
$$
Here, according to Proposition 2.2, the first term on the right-hand side 
is minimized and is equal to 1 when $U = \frac{1}{2}\,Z$ where $Z$ 
has density $f_\alpha$. In that case, for the second term we have
$$
\frac{A_\alpha}{N_\alpha(U)} = \frac{1}{\Var(U)} = \frac{4}{\Var(Z)} = 
\frac{4(3\alpha - 1)}{\alpha - 1},
$$
where we recalled the identity (2.6). Thus, we arrive at the following
relation, which contains (3.2) for $\alpha = \infty$.

\vskip5mm
{\bf Proposition 3.3.} {\sl Let $1 < \alpha \leq \infty$. For any
integer valued random variable $X$ having finite variance, we have
\be
1 \, \leq \, N_\alpha(X) \, \leq \,
1 + \frac{4(3\alpha - 1)}{\alpha - 1}\,\Var(X).
\en
}

\vskip5mm
\section{{\bf Log-concave Sequences}}
\setcounter{equation}{0}

\vskip2mm
\noindent
We now turn to lower bounds for the R\'enyi entropies. This cannot be performed 
in terms of variances in the entire class of discrete probability distributions, 
so some extra hypotheses are needed. As it turns out,
the class of discrete log-concave distributions perfectly fits our aims.
First, let us recall that a sequence $\{a_k\}_{k \in \Z}$ of non-negative
numbers is called log-concave, if
\be
a_k^2 \geq a_{k-1}\, a_{k+1}
\en
for all $k \in \Z$. Similarly, a finite sequence $\{a_k\}_{k = m}^n$ of 
non-negative numbers is log-concave, if this inequality is fulfilled
whenever $m+1 \leq k \leq n-1$. Defining $a_k = 0$ for $k<m$ and $k>n$,
we then obtain an infinite log-concave sequence.

Log-concave sequences appear in a purely algebraic framework.
The following classical result goes back to Newton. 

\vskip5mm
{\bf Proposition 4.1.}
{\sl Suppose that a polynomial $P(z) = \sum_{k=0}^n {n \choose k}\, a_k z^k$ 
has only real zeros over the field $\C$ of complex numbers. 
Then $(4.1)$ holds true for all $1 \leq k \leq n-1$.
}

\vskip5mm
Following \cite{S}, one may give a proof based
on the following observation of independent interest due to Gauss-Lucas:
For any polynomial $P$ over $\C$ of degree ${\rm deg}(P) \geq 1$, the zeros 
of its derivative $P'$ belong to the convex hull of the zeros of $P$.

The sequences that produce polynomials with real roots give an interesting 
subset of the class of sequences satisfying a log-concavity inequality. 
One should note that such sequences need not have contiguous support. 
For example, $z^2 -1$ has roots $\pm 1$ and corresponds to the sequence $\{1,0,-1\}$.

Let us now mention the following important theorem due to Hoggar \cite{H}.

\vskip5mm
{\bf Proposition 4.2.} {\sl If the coefficients of two polynomials are positive 
and form log-concave sequences, then so does their product.
}

\vskip5mm
Note that for
$P(z) = \sum_{k=0}^n a_k z^k$ and $Q(z) = \sum_{k=0}^m b_k z^k$ of degrees
$n$ and $m$ respectively,
$$
R(z) \, = \, P(z) Q(z) \, = \, \sum_{k=0}^{n+m} c_k z^k, \quad
c_k \ = \sum_{k_1 + k_2 = k} a_{k_1} b_{k_2}.
$$
Here, the coefficients $\{c_k\}_{k=0}^{n+m}$ appear as the convolution of the 
sequences $\{a_k\}_{k = 0}^n$ and $\{b_k\}_{k = 0}^m$. Thus, Proposition 4.2
tells us that the class of finite positive log-concave sequences
is closed under the convolution operation.

The assumption about the positivity may not be removed in this conclusion.
For a counter-example, one may take the polynomials
$$
P(z) = 1 + z, \quad Q(z) = 1+z^3, \quad R(z) = 1 + z + z^3 + z^4,
$$
which correspond to the sequences $\{1,1\}$, $\{1,0,0,1\}$, $\{1,1,0,1\}$.
Here, the first two sequences are log-concave, while the third one is not. 
With this in mind, the notion of a log-concave discrete distribution
on $\Z$ should be introduced as in Introduction.

Recalling the definition (1.3),
as an immediate consequence from Proposition 4.2, we obtain:

\vskip5mm
{\bf Corollary 4.3.} {\sl The class of (discrete) log-concave probability 
distributions on $\Z$ is closed under the convolution operation.}

\vskip5mm
\section{{\bf From Discrete to Continuous Log-concave Measures}}
\setcounter{equation}{0}

\vskip2mm
\noindent
Let us turn to basic properties of discrete log-concave distributions.
Some of them can be obtained by employing known results from the theory of 
``continuous'' log-concave functions. To this aim, we describe a simple
construction, which allows one to associate with a discrete log-concave
distribution defined by a probability function
a certain log-concave function on the real line. 

Given a random variable $X$ with a discrete log-concave distributions defined by a probability function $f(k)$, denote by $\Delta$ the smallest closed interval
containing ${\rm supp}(f)$. Let us extend the function
$$
V(k) = -\log f(k)
$$
linearly on every segment $[k,k+1] \subset \Delta$, and put $V = \infty$ 
outside $\Delta$ (if $\Delta$ is not the whole real line). Then $V$ is finite 
and convex on $\Delta$. To see this, assume that 
$|\Delta| \geq 2$, and for $[k-1,k+1] \subset \Delta$, rewrite (4.1) as
$$
V(k) - V(k-1) \leq V(k+1) - V(k).
$$
But $V(k) - V(k-1) = V'(x)$ on $(k-1,k)$, which means that $V$ has a non-decreasing
(Radon-Nikodym) derivative. As a result, we obtain a log-concave function
$$
f(x) = e^{-V(x)}, \quad x \in \R,
$$
which coincides with $f(k)$ on all integers $k$. Let us call it
a log-piecewise linear extension of the sequence $f(k)$. 

In general, $f(x)$ does not need be a probability density on the line. Nevertheless, 
since
$$
\int_k^{k+1} f(x)\,dx \leq \max\{f(k),f(k+1)\} \leq f(k) + f(k+1),
$$
after summation over all $k \in \R$, we conclude that $f$ is integrable. 
Let us emphasize this fact once more.

\vskip5mm
{\bf Proposition 5.1.} {\sl The restrictions of densities of finite (continuous)
log-concave measures from $\R$ to $\Z$ describe the whole class of
probability functions of discrete log-concave measures.
}

\vskip5mm
Returning to the log-piecewise linear extension of the sequence $f(k)$
and using an additional property $f(k) \leq 1$, an immediate consequence 
of the integrability of $f$ is that, with some constant $c>0$ depending on 
the distribution of $X$, 
$$
f(k) \leq 2\,e^{-c |k|} \quad {\rm for \ all} \ \ k \in \Z.
$$
In particular, there is a point of maximum for this sequence, which is called a mode
(using a probabilistic language). It is also a mode for the function $f(x)$.
So, it is a mode of $X$.

Let $m \in \Z$ be a mode of $X$. Then necessarily $f(k)$ is non-increasing in 
$k \geq m$ and is non-decreasing in $k \leq m$. Hence
\bee
\int_k^{k+1} f(x)\,dx \, \geq \, f(k+1) \ \ {\rm for} \ k \geq m, \qquad
\int_{k-1}^k f(x)\,dx \, \geq \, f(k-1) \ \ {\rm for} \ k \leq m.
\ene
Performing summation over all $k$, we get
$$
\int_{-\infty}^\infty f(x)\,dx \geq 1 - f(m).
$$
Also, since $V$ is linear on each interval $[k,k+1] \subset \Delta$, 
$f(x)$ is convex, and therefore its graph lies below the segment
connecting the points $(k,f(k))$ and $(k+1,f(k+1))$. Hence
$$
\int_k^{k+1} f(x)\,dx \leq \frac{f(k) + f(k+1)}{2}.
$$
Performing summation over all $k$, we arrive at:

\vskip5mm
{\bf Proposition 5.2.} {\sl If $m$ is mode of $X$ and $|\Delta| \geq 2$, then
$$
1 - f(m) \leq \int_{-\infty}^\infty f(x)\,dx \leq 1.
$$
}

\vskip2mm
Similar bounds hold true for the second moment. In particular,
one may easily derive:

\vskip5mm
{\bf Proposition 5.3.} {\sl Let $X$ be a discrete random variable with 
a log-concave density $f$ with mode at $m$. If $f(m) \leq \frac{1}{2}$, then
the log-piecewise extension of the density satisfies
$$
\frac{1}{3}\, \E\, (X-m)^2 \leq \int_{-\infty}^\infty (x-m)^2\, f(x)\,dx \leq  
3\, \E\, (X-m)^2.
$$
Moreover, the left inequality is also true when $f(m) \geq \frac{1}{2}$.
}

\vskip5mm
\section{{\bf Second Moment and Maximum of Density}}
\setcounter{equation}{0}

\vskip2mm
\noindent
Propositions 5.2-5.3 may be used to extend a number of results about continuous 
log-concave  distributions to the discrete setting. Let us remind basic relations 
for the class of log-concave densities on the line in the usual continuous setting.
One important feature in this case is that the general lower bound 
as in Proposition 2.1 may be reversed.

Keeping $X$ to denote a discrete random variable, let $Z$ be a (continuous) 
random variable with a log-concave density $g$. It is supported on some 
closed interval $\Delta \subset \R$, finite or not, which contains a point 
of maximum $m$ of $g$, i.e., the mode of $Z$ (which may be one of the endpoints 
of $\Delta$ when this interval is bounded from the left or from the right). 
Hence $M(Z) = g(m)$. Together with (2.2), it is known that
\be
\frac{1}{12} \leq M^2(Z)\,\Var(Z) \leq 1,
\en
where an inequality on the right-hand side is achieved for 
the one-sided exponential distribution with density 
$g(x) = e^{-x}\,1_{(0,\infty)}(x)$ (cf. Proposition 2.1 in \cite{B-G3}). 
In fact, if $g$ is symmetric about $m$, the upper bound can be improved to
\be
M^2(Z)\,\Var(Z) \leq \frac{1}{2},
\en
which is attained for the two-sided exponential distribution with density 
$g(x) = \frac{1}{2}\,e^{-|x|}$.

There is a similar assertion about 
the second moment when $Z$ is centered at the mode. 

\vskip5mm
{\bf Proposition 6.1.} {\sl If a (continuous) random variable $Z$
has a log-concave density $g$ with mode at the point $m$, then
$$
M^2(Z)\,\E\, (Z-m)^2 \leq 2.
$$
}

\vskip2mm
Combining this relation with the lower bound in (6.1), we also get
$$
\E\, (Z-m)^2 \leq 24\,\Var(Z).
$$

Let us now describe an immediate application of Proposition 6.1
to the discrete setting. Suppose that we are given an integer-valued random 
variable $X$ with discrete log-concave distribution, and denote by $f(x)$ 
the log-piecewise linear extension of the probability function 
$f(k) = \P\{X=k\}$, $k \in \Z$. Let $Z$ be a random variable 
whose density $g$ is proportional to $f$, that is, $g(x) = \frac{1}{B}\,f(x)$, 
$B = \int_{-\infty}^\infty f(x)\,dx$. Recall that, by Proposition 5.2, 
$B \leq 1$. Hence, by Propositions 6.1 applied to $g$, we get that
$$
f(m)^2 \int_{-\infty}^\infty (x-m)^2\, f(x)\,dx \leq 2B^3.
$$
On the other hand, according to Proposition 5.3, the above integral
dominates $\frac{1}{3}\, \E\, (X-m)^2$.
Combining the two inequalities, we arrive at:

\vskip5mm
{\bf Proposition 6.2.} {\sl If $X$ is a discrete random variable with 
a log-concave density $f$ with mode at $m$, then
$$
M^2(X)\,\E\, (X-m)^2 \leq 6.
$$
In particular,
\be
M^2(X)\,\Var(X) \leq 6.
\en
}

\vskip2mm
{\bf Proof of Proposition 6.1.} 
We may assume that $m=0$ and $g(0)=1$. Moreover, let
$a = \P\{Z < 0\}$ and $b = \P\{Z>0\}$ be positive.

By the log-concavity, $g(x) = \exp\{-V(x)\}$ for some convex function
$V \colon \R \rightarrow [0,\infty]$ with $V(0) = 0$. Putting
$\tilde g(x) = \exp\{-\tilde V(x)\}$ with
$\tilde V(x) = x/b$ for $x \geq 0$ and $\tilde V(x) = -x/a$ for $x \leq 0$,
we obtain another log-concave probability density such that
$$
\int_{-\infty}^0 g(x)\,dx = \int_{-\infty}^0 \tilde g(x)\,dx = a, \qquad
\int_0^\infty g(x)\,dx = \int_0^\infty \tilde g(x)\,dx = b.
$$
Since $\tilde V$ is linear on $[0,\infty)$ and $\tilde V(0) = 0$,
necessarily $V(x) \leq \tilde V(x)$ on the interval $0 \leq x < x_0$
and $V(x) \geq \tilde V(x)$ on the half-axis $x > x_0$ for some $x_0 > 0$.
Equivalently, $g(x) \geq \tilde g(x)$ for $0 \leq x < x_0$ and
$g(x) \leq \tilde g(x)$ for $x > x_0$. Introduce
the distribution functions 
$$
G(x) = \P\{Z \leq x\} = \int_{-\infty}^x g(y)\,dy \quad {\rm and} \quad
\tilde G(x) = \int_{-\infty}^x \tilde g(y)\,dy.
$$
It follows that the function $\psi(x) = G(x) - \tilde G(x)$ is vanishing
at the origin and at infinity, and its derivative is non-negative 
on $[0,x_0)$ and non-positive on $(x_0,\infty)$. Hence $\psi(x) \geq 0$
for all $x \geq 0$. This implies that
$$
\int_0^\infty x^2\,g(x)\,dx - \int_0^\infty x^2\,\tilde g(x)\,dx
= -2 \int_0^\infty x \psi(x)\,dx \leq 0.
$$
By the same argument, a similar inequality holds when integrating
over the negative half-axis. Hence
$$
\E\, Z^2 \, = \, \int_{-\infty}^\infty x^2\,g(x)\,dx \, \leq \,
\int_{-\infty}^\infty x^2\,\tilde g(x)\,dx \, = \, 2\,(a^3 + b^3) \, \leq \, 
2\,(a + b)^3 \, = \, 2.
$$
\qed

\vskip5mm
{\bf Remark.} In terms of the R\'enyi entropy powers, the inequalities
(6.1)-(6.2) take the form
\be
\frac{1}{12}\, N_\infty(Z) \leq \Var(Z) \leq N_\infty(Z), \qquad 
\Var(Z) \leq \frac{1}{2}\, N_\infty(Z).
\en
Here, the first inequality is general and corresponds to (2.4)
in the limit case $\alpha = \infty$, while the upper bounds on the variance hold
for log-concave probability measures (symmetric in the second case).
Similarly to (2.4), these upper bounds may be extended to the whole
range of $\alpha$ at the expense of certain $\alpha$-dependent factors.
Although the extremal log-concave distributions are not known so far,
one particular case may easily be settled. 
Given a random variable $Z$ with a log-concave distribution,
let us apply the second upper bound in (6.4) to $\tilde Z = Z - Z'$,
where $Z'$ is an independent copy of $Z$. Since
$N_\infty(\tilde Z) = N_2(Z)$, we then obtain that
$$
\Var(Z) \leq \frac{1}{4}\, N_2(Z).
$$
Here, we have an equality for the exponential distribution, and thus
 the conclusion that among continuous log-concave variables on the real
 line, the exponential distribution has minimal $2$-R\'enyi entropy, often
 referred to as collision entropy.

\vskip5mm
\section{{\bf Sharpenings via Rearrangement Arguments}}
\setcounter{equation}{0}

\vskip2mm
\noindent
The inequality (6.3) may be sharpened and generalized by involving more
delicate arguments based on the so-called rearrangement of densities.
Similarly to the continuous setting, these arguments allow one 
to explore an extremal role of a discrete counterpart of the exponential
distributions. More precisely, when a probability mass function 
$f \colon \Z \to [0,\infty)$ can be written as $f(k) = C p^{|k|}$ for some
$0 \leq p < 1$, we will call it a symmetric two-sided geometric distribution.
Here, the symmetry property refers to the identity $f(-k) = f(k)$. 

Let us start with basic definitions.

\vskip5mm
{\bf Definition} (Decreasing rearrangement).
For a probability function $f \colon \Z \to [0,\infty)$, denote by 
$f^{\downarrow}$ its decreasing rearrangement. Explicitly, $f^{\downarrow}$ 
is defined on $\N$ and satisfies $f^{\downarrow}(k) \geq f^{\downarrow}(k+1)$ 
for all $k \geq 0$ and there exists a bijection $\tau \colon \N \to \Z$ such that 
$f^{\downarrow}(k) = f(\tau(k))$.

\vskip5mm
{\bf Definition} (Majorization).
Given probability functions $f$ and $g$ on $\Z$, we say that $f$ majorizes 
$g$ and write $f \succ g$, when
$\sum_{k=0}^n f^{\downarrow}(k) \geq  \sum_{k=0}^n g^{\downarrow}(k)$
for all $n \geq 0$.

\vskip5mm
We extend this notion to random variables by writing $X \succ Y$ when 
$f \succ g$ holds for their respective probability functions.

\vskip5mm
{\bf Definition} (Schur-Concavity).
A functional $\Phi$ defined on a given family of probability functions on $\N$
is called Schur-convex, when $f \succ g$ implies $\Phi(f) \geq \Phi(g)$. 
A functional $\Phi$ is Schur-concave when $-\Phi$ is Schur-convex.

\vskip5mm
We extend Schur-convexity/concavity to random variables by asking that 
Schur-convexity/ concavity hold for the respective probability mass functions.

\vskip5mm
{\bf Lemma 7.1} (\cite{M-T2}). {\sl Any symmetric log-concave probability function 
$f$ majorizes the symmetric two-sided geometric distribution $g$ with 
the same maximum as $f$.
}

\vskip5mm
A more general result was actually first proven in \cite{M-T}. We include a proof of 
Lemma 7.1 for completeness. 

\vskip2mm
{\bf Proof.}
Let $q \in [0,1]$. Define, for $k \in \Z$,
$
g_q(k) = q^{|k|} f(0)
$
with the convention that $0^0 = 1$. Note that
$$ 
\sum_{k \in \Z} g_0(k) = f(0) \leq \sum_{k \in \Z} f(k) = 1 
$$
and $\sum_{k \in \Z} g_1(n) = \infty$. Since the function 
$q \mapsto \sum_{k \in \Z} g_q(k)$ is continuous and strictly increasing on 
$[0,1]$, we deduce by the intermediate value theorem that there exists 
a unique $q^* \in [0,1)$ such that
$$ 
\sum_{k \in \Z} g_{q^*}(k) = \sum_{k \in \Z} f(k) = 1. 
$$
By construction, $g_{q^*}$ is the probability function of a two-sided 
geometric distribution. Denote by the same letters $f$ and $g_q$ the corresponding
log-piecewise linear extensions. In particular,
$$ 
\max_{x \in \R} g_{q^*}(x) = g_{q^*}(0) = f(0) = \max_{x \in \R} f(x). 
$$
Thus, the function $L = -\log g_{q^*}$ is affine on $[0,\infty)$, 
while $V = -\log f$ is convex, with $L(0)=V(0)$. Hence,
there must exist $l \geq 0$ such that $V \leq L$ on $[0,l]$ and 
$V \geq L$ on $[l,\infty)$, that is, $f \geq g_{q^*}$ on $[0,l]$ and 
$f \leq g_{q^*}$ on $[l,\infty)$. Hence, for all $n \geq 0$,
$\sum_{k=0}^n f^{\downarrow}(k) \geq \sum_{k=0}^n g_{q^*}^{\downarrow}(k)$.
\qed

\vskip5mm
{\bf Lemma 7.2} (Schur-concavity of moments on symmetric log-concaves).
{\sl If the function $u \colon [0,\infty) \rightarrow [0,\infty)$ is non-decreasing, 
then the functional $\Phi(f) = \E\,u(|X|)$ is Schur concave on the set of 
symmetric discrete log-concave probability functions.
}

\vskip5mm
{\bf Proof}. By linearity of $\Phi$ with respect to $u$,
one may assume that $u(x) = 1_{(\lambda,\infty)}(x)$ for some integer 
$\lambda \geq 0$. Suppose that $X$ and $Y$ are symmetric log-concave variables 
such that $f_X \succ f_Y$. By symmetry, unimodality and the assumed majorization, 
we have
\bee
\E\,u(|X|) 
 & = & \P\{|X| > \lambda\} \ = \
1 - \P\{|X| \leq \lambda\} \ = \
1 - \sum_{k = -\lambda}^{\lambda} f_X(k) \\
 & = &
1 - \sum_{k=0}^{2\lambda} f_X^{\downarrow}(k)
 \ \leq \
1 - \sum_{k=0}^{2\lambda} f_Y^{\downarrow}(k)
 \ = \
\P\{|Y| > \lambda\} \ = \ \E\,u(|Y|),
\ene
and the result follows.
\qed

\vskip5mm
We are prepared to settle the upper bound of Theorem 1.1 in the
symmetric case (and in somewhat more general form).

\vskip5mm
{\bf Proposition 7.3.} {\sl Let $\Psi \colon [0,\infty) \to [0,\infty)$ be a 
non-decreasing function. For any symmetric discrete log-concave random variable $X$,
$$ 
M(X)\, \Psi(\Var(X)) \, \leq \, \sup_{0 \leq p < 1} \, 
\bigg[\frac{1-p}{1+p}\ \Psi\Big(\frac{2p}{(1-p)^2}\Big)\bigg]. 
$$
In particular,
\be
M^2(X)\,(1 + 2\,\Var(X)) \leq 1. 
\en
The latter inequality becomes an equality for all two-sided 
geometric distributions. 
}

\vskip5mm
{\bf Proof.}
Given $X$ a symmetric log-concave random variable, let $Y$ be the two sided 
geometric distribution with same maximum (necessarily attained at 0). 
Then by Lemma 7.1, we have $X \succ Y$, so that by Lemma 7.2, 
$\Var(X) \leq \Var(Y)$. Thus
$$ 
M(X)\, \Psi(\Var(X)) \, \leq \, M(Y)\, \Psi(\Var(Y)). 
$$
Note that for a two sided geometric distribution with probability mass 
function $g(k) = Cp^{|k|}$, $0 \leq p < 1$, one has
$$ 
C^{-1} = \sum_{k \in \Z} p^{|k|} = \frac{1+p}{1-p}, \qquad
M(Y) = C = \frac{1-p}{1+p}, 
$$
$$ 
\Var(Y) = \E Y^2 = \frac{2p}{(1-p)^2}. 
$$
Hence
$$ 
M(Y)\, \Psi(\Var(Y)) \, \leq \,
\sup_{p \in [0,1)} \bigg[\frac{1-p}{1+p}\ \Psi\Big(\frac{2p}{(1-p)^2}\Big)\bigg]. 
$$

In particular, for $\Psi(x) = \sqrt{1 + 2x}$, the function under the above 
supremum is equal to 1 for any value of the parameter $p$, and we obtain (7.1)
together with the assertion about the extremal role of two-sided geometric
distributions.
\qed

\vskip5mm
While (7.1) is a full analog of the inequality (6.2) from the continuous
setting under the symmetry hypothesis, the more general inequality (6.1) also
extends to the discrete setting. The next assertion provides the upper bound 
of Theorem 1.1 in the non-symmetric case.

\vskip5mm
{\bf Proposition 7.4.} {\sl For any discrete log-concave random variable $X$,
\be
M^2(X)\,(1 + 4\,\Var(X)) \leq 4. 
\en
This bound is asymptotically attained for one-sided
geometric distributions. In particular,
\be
M^2(X)\,\Var(X) \leq 1. 
\en
}

\vskip2mm
{\bf Proof.} We employ another relation,
\be
H_\alpha(X) \leq H_\infty(X) + \log\big(\alpha^{\frac{1}{\alpha - 1}}\big),
\qquad 0 < \alpha < \infty,
\en
which was recently obtained in \cite{M-T2} using the rearrangement arguments. 
This inequality holds for any discrete log-concave random variable $X$ and 
is asymptotically attained for the geometric distributions with probability 
functions $f(k) = (1-p)\,p^k$, $k \geq 0$, as $p \rightarrow 1$.
In particular, for $\alpha=2$, it takes the form
$H_2(X) \leq H_\infty(X) + \log 2$, which is the same as
\be
M(X) \, \leq \, 2\, \sum_{k \in \Z} f(k)^2
\en
in terms of the probability function $f(k)$ of $X$.

Note that, if $Y$ is an independent copy of $X$, the random variable $X-Y$
will have a symmetric discrete log-concave distribution with
\be
M(X-Y) = \P\{X-Y = 0\} = \sum_{k \in \Z} f(k)^2,
\en
so that, by (7.5),
\be
M^2(X) \, \leq \, 4\, M^2(X-Y).
\en

On the other hand, by Proposition 7.3, 
$M^2(X-Y)\,(1 + 2\, \Var(X-Y)) \leq 1$, that is,
$$
M^2(X-Y)\,(1 + 4\, \Var(X)) \leq 1.
$$
It remains to apply this inequality in (7.7).
\qed

\vskip2mm
Note that for the geometric distribution with probability function 
$f(k) = (1-p)\,p^k$, $k \geq 0$, we have $M^2(X)\,(1 + 4\,\Var(X)) = (1+p)^2$. 
Hence, (7.2) is asymptotically attained for $p \rightarrow 1$.

\vskip5mm
\section{{\bf Concentration Functions}}
\setcounter{equation}{0}

\vskip2mm
\noindent
Turning to applications,
first let us relate the concentration function to the $M$-functional.

\vskip5mm
{\bf Lemma 8.1.} {\sl For any random variable $X$,
\be
Q(X;\lambda) \, = \, \lambda M(X+U_\lambda), \quad \lambda > 0,
\en
where the random variable $U_\lambda$ is independent of $X$ and has a uniform
distribution on the interval $(0,\lambda)$. Analogously, in the discrete
setting, we have
\be
Q(X;\lambda) \, = \, (\lambda + 1)\, M(X+U_\lambda), \quad 
\lambda = 0,1,2,\dots,
\en
assuming that the random variable $U_\lambda$ is independent of $X$ and has 
a discrete uniform distribution on the integer interval $\{0,1,\dots,\lambda\}$.
}

\vskip5mm
{\bf Proof}. The first claim follows from the fact that the random variable 
$X_\lambda = X+U_\lambda$ has an absolutely continuous distribution with density
$$
f_\lambda(x) = \frac{1}{\lambda}\,\P\{x-\lambda \leq X \leq x\} \quad {\rm a.e.}
$$
According to the definition (1.1), this formula yields
$$
M(X_\lambda) \, = \, {\rm ess\,esup}_x\, f_\lambda(x) \, = \,
\frac{1}{\lambda}\,Q(X;\lambda).
$$

For the second claim, we similarly have that $X_\lambda$ takes integer values
with probabilities
$$
f_\lambda(k) = \frac{1}{\lambda + 1}\,\P\{k-\lambda \leq X \leq k\}, \quad k \in \Z.
$$
Since the supremum in (1.1) is attained for some integer value $x = k$,
Lemma 8.1 follows.
\qed

\vskip5mm
One can now apply the lower bound of Proposition 2.1, cf. (2.2), to the 
random variables $X_\lambda$ in (8.1). This leads to a corresponding lower bound 
for the concentration function (which is actually known, cf. \cite{B-G3}).

\vskip5mm
{\bf Proposition 8.2.} {\sl For any random variable $X$,
\be
Q(X;\lambda) \, \geq \, \frac{\lambda}{\sqrt{\lambda^2 + 12\,\Var(X)}}, 
\quad \lambda > 0.
\en
}

\vskip2mm
As a consequence, applying (8.3) with $\lambda \uparrow 1$
in the discrete setting, we arrive at the lower bound (1.4) in Theorem 1.1.
Note that since in general $M(X) \leq 1$, while the variance may take any prescribed 
value (already within specific families such as Poisson distributions),
the constant $1$ may not be removed from the square root in (1.4). 

More generally, if the random variable $X$ is discrete and $\lambda \geq 0$ 
is integer, we apply (8.3) with $\lambda' \uparrow \lambda+1$, 
so that to get in the limit
\be
Q(X;\lambda) \, \geq \, \frac{\lambda + 1}{\sqrt{(\lambda + 1)^2 + 12\,\Var(X)}}.
\en

For the upper bound, one may start from the identity (8.2) and apply 
the upper bounds of Theorem 1.1 to the random variables $X_\lambda = X+U_\lambda$ 
in place of $X$, assuming that the random variable $U_\lambda$ is independent 
of $X$ and has a discrete uniform distribution on the integer interval 
$\{0,1,\dots,\lambda\}$. Note that $X_\lambda$ is symmetric around the point
$\lambda/2$, if $X$ is symmetric about the origin. Since
$$
\Var(U_\lambda) = \frac{\lambda (\lambda + 2)}{12},
$$
together with (8.4) we are lead to the following statement.

\vskip5mm
{\bf Proposition 8.3.} {\sl If the random variable $X$ has a log-concave discrete 
distribution, then, for any integer $\lambda \geq 0$,
$$
\frac{\lambda + 1}{\sqrt{1 + \lambda (\lambda + 2) + 12\,\Var(X)}} 
 \, \leq \, Q(X;\lambda) \, \leq \, 
\frac{2\,(\lambda+1)}{\sqrt{1 + \frac{\lambda (\lambda + 2)}{3} + 4\,\Var(X)}}.
$$
Moreover, if the distribution of $X$ is symmetric about a point, then the upper
bound may be sharpened to
$$
Q(X;\lambda) \leq 
\frac{\lambda + 1}{\sqrt{1 + \frac{\lambda (\lambda + 2)}{6} + 2\,\Var(X)}}.
$$
}

\vskip5mm
Here, the value $\lambda = 0$ returns us to the statement of Theorem 1.1.

\vskip5mm
\section{{\bf Proof of Theorems 1.2-1.3}}
\setcounter{equation}{0}

\vskip2mm
\noindent
Let $X_1,\dots,X_n$ be independent discrete random variables, and
$S_n = X_1 + \dots + X_n$.

\vskip5mm
{\bf Proof of Theorem 1.2.} First, consider the symmetric case.
By Proposition 3.3, for any integer valued 
random variable $X$ having a finite variance,
$$
\Delta_\alpha(X) \, \leq \, \frac{4(3\alpha - 1)}{\alpha - 1}\,\Var(X).
$$
On the other hand, if the distribution of $X$ is symmetric about a point
and log-concave, the inequality (1.5) of Theorem 1.1 yields a lower bound
\be
\Delta_\infty(X) \geq 2\, \Var(X).
\en
Since the function $\alpha \rightarrow \Delta_\alpha(X)$ is non-increasing, 
the two inequalities give
\be
2\,\Var(X) \, \leq \, 
\Delta_\alpha(X) \, \leq \, \frac{4(3\alpha - 1)}{\alpha - 1}\,\Var(X).
\en
Applying this to $X_k$'s, as well as to $X = S_n$, we then get
$$
\Delta_\alpha(S_n) \, \geq \, 2\,\Var(S_n) \, = \, 2 \, \sum_{k=1}^n \Var(X_k)
\, \geq \, 
\frac{\alpha - 1}{2(3\alpha - 1)}\ \sum_{k=1}^n \Delta_\alpha(X_k).
$$

In fact, up to an $\alpha$-dependent constant, this upper bound may be
reversed, since, by (9.2), we also have
\bee
\Delta_\alpha(S_n) 
 & \leq &
\frac{4(3\alpha - 1)}{\alpha - 1}\ \Var(S_n) \\
 & \leq &
\frac{4(3\alpha - 1)}{\alpha - 1}\ \sum_{k=1}^n \Var(X_k) \, \leq \,
\frac{2(3\alpha - 1)}{\alpha - 1}\ \sum_{k=1}^n \Delta_\alpha(X_k).
\ene
This proves the desired relations in (1.10).

Now, suppose that $1 < \alpha \leq 2$ and drop the symmetry assumption.
If $Y$ is an independent copy of $X$, the random variable $\tilde X = X-Y$ 
will have a symmetric discrete log-concave distribution. Hence, (9.1)
is applicable to $\tilde X$. On the other hand, the identity (7.6) reads as 
$N_\infty(\tilde X) = N_2(X)$, that is, $\Delta_\infty(\tilde X) = \Delta_2(X)$.
Hence, from (9.1) we get that
$$
\Delta_\alpha(X) \geq \Delta_2(X) \geq 4\, \Var(X).
$$
Combining this with the upper bound in (9.2), we arrive at its sharpened
variant
\be
4\,\Var(X) \, \leq \, 
\Delta_\alpha(X) \, \leq \, \frac{4(3\alpha - 1)}{\alpha - 1}\,\Var(X).
\en

It remains to apply this to $X_k$'s, as well as to $X = S_n$ similarly as above,
and then we obtain the lower bound as in (1.10) with an improved constant.
The proof of the upper bound is based on (9.3) and is also similar.
\qed

\vskip5mm
{\bf Remark.}
In the case $\alpha = \infty$, (1.10) takes the form
$$
\frac{1}{6}\,\sum_{k=1}^n \Delta_\infty(X_k) \, \leq \, 
\Delta_\infty(S_n) \, \leq \, 6\,\sum_{k=1}^n \Delta_\infty(X_k).
$$
Here, the lower bound reminds the relation
$$
N_\infty(S_n) \, \geq \, \frac{1}{2} \sum_{k=1}^n N_\infty(X_k)
$$
for the continuous setting (without any assumptions on the shape of
distributions, cf. \cite{B-G1,M-M-X}). Since the functions
$\alpha \rightarrow \Delta_\alpha(X)$ are non-increasing, it implies that
$$
\Delta_1(S_n) \, \geq \, \frac{1}{6} \sum_{k=1}^n \Delta_\infty(X_k).
$$
However, one cannot estimate $\Delta_\infty(X_k)$ in terms of $\Delta_1(X_k)$.
This can be seen on the example of the two-sided geometric distribution with
probability function $f(k) = C p^{|k|}$, $k \in \Z$, for small $0<p<1$. 
In this case, 
$$
\Delta_\infty(X) = 2\,\Var(X) = \frac{4p}{(1-p)^2},
$$
which corresponds to the equality in (1.5), while
$$
H(X) = -\log(1-p) + \log(1+p) - \frac{2p \log p}{(1-p)(1+p)}.
$$
As $p \rightarrow 0$, we have $H(p) \sim 2p\, \log(1/p)$, hence
\bee
\Delta(X) \, = \, N(X) - 1 
 & = & 
e^{2H(X)} - 1 \\
 & \sim & 
2H(X) \, \sim \, 4p\, \log(1/p) \, \sim \, 2\,\Var(X)\, \log(1/p). 
\ene

\vskip5mm
{\bf Proof of Theorem 1.3.} Suppose that $\Var(X_k) \geq \sigma^2$ 
for all $k \leq n$ with some $\sigma>0$. 

By Proposition 3.1, if $X$ is a discrete
random variable with $\Var(X) \geq \sigma^2$, we have
$$
N_\alpha(X) \, \leq \, 2\pi e\,\Big(1 + \frac{1}{12 \sigma^2}\Big)\, \Var(X)
$$
for any $\alpha \geq 1$. In addition, according to the upper bound (1.4)
of Theorem 1.1,
\be
N_\infty(X) \geq \frac{1}{4} + \Var(X) \geq \Var(X).
\en
Since $\alpha \rightarrow N_\alpha(X)$ is a non-increasing function,
these bounds imply that
\be
\Var(X) \, \leq \, N_\infty(X) \, \leq \, N_\alpha(X) \, \leq \,
2\pi e\,\Big(1 + \frac{1}{12 \sigma^2}\Big)\, \Var(X).
\en
Being applied first to $X=S_n$ and then to each $X = X_k$, (9.5) yields
$$
N_\alpha(S_n) \, \geq \, \Var(S_n) \, = \, \sum_{k=1}^n \Var(X_k) \, \geq \,
\frac{1}{2\pi e\,(1 + \frac{1}{12 \sigma^2})}\,\sum_{k=1}^n N_\alpha(X_k).
$$
This proves the first assertion (1.11) of the theorem.

By a similar argument, the resulting bound may be reversed, even without
any condition on variances, and actually further strengthened.
Indeed, applying (9.5) to $X =S_n$ and using the first inequality in (9.4)
in the weaker form $N_\alpha(X_k) \geq \frac{1}{4} + \Var(X_k)$, we have
\bee
N_\alpha(S_n) 
 & \leq &
2\pi e\,\Big(\frac{1}{12} + \Var(S_n)\Big) \\
 & = &
2\pi e\,\Big(\frac{1}{12} + \sum_{k=1}^n \Var(X_k)\Big) \ \leq \
2\pi e\,\Big(\frac{1}{12} + \sum_{k=1}^n \Big(N_\alpha(X_k) - \frac{1}{4}\Big)\Big),
\ene
and (1.12) immediately follows.
\qed

\vskip5mm
\section{{\bf Remarks on Bernoulli sums}}
\setcounter{equation}{0}

\vskip2mm
\noindent
Finally, let us illustrate Theorem 1.1 on the example of the Bernoulli sums
$$
S_n = X_1 + \dots + X_n,
$$
where the independent summands $X_k$ take the values 1 and 0 with probabilities 
$p_k$ and $q_k$ respectively. By Hoggar's theorem on preservation of 
log-concavity under convolutions, the distribution of $S_n$ is discrete 
log-concave (Corollary 4.3). This property is not so obvious 
on the basis of the explicit expression for the probability function 
$$
f_n(k) = \P\{S_n = k\} = 
\sum p_1^{\ep_1} q_1^{1 - \ep_1} \dots p_n^{\ep_1} q_n^{1 - \ep_n}, \qquad
k = 0,1,\dots,n,
$$
where the summation is running over all $0-1$ sequences $\ep_1,\dots,\ep_n$
such that $\ep_1 + \dots + \ep_n = k$.

One of the challenging problems about $S_n$ has been to how effectively estimate 
from above the maximum $M(S_n) = \max_k f_n(k)$ in terms of $p_k$'s. 
In this particular model, there exist several different approaches to the
problem. First let us describe the standard approach in Probability Theory which 
is irrelevant to log-concavity. In general, the concentration function of 
a random variable $X$ may be bounded in terms of the characteristic function 
$v(t) = \E\,e^{itX}$ by virtue of Esseen's bound
\be
Q(X;\lambda) \, \leq \, 
\Big(\frac{96}{95}\Big)^2 \, \lambda \int_{-1/\lambda}^{1/\lambda} |v(t)|\,dt,
\qquad \lambda>0,
\en
cf. e.g. \cite{P1}, \cite{P2}. For the characteristic function of $S_n$ we have
$v(t) = \prod_{k=1}^n (q_k + p_k e^{it})$. Since 
$$
|q + p e^{it}|^2 = p^2 + q^2 + 2pq\,\cos t = 1 - 4pq\,\sin^2(t/2) \leq 
\exp\{- 4pq\,\sin^2(t/2)\},
$$
it follows that
$$
|v(t)| \leq \exp\Big\{- 2\sum_{k=1}^n p_k q_k\,\sin^2(t/2)\Big\} = 
\exp\Big\{- 2\sigma^2 \sin^2(t/2)\Big\},
$$
where $\sigma^2 = \Var(S_n)$ ($\sigma > 0$).
Using $|\sin x| \geq \frac{2}{\pi}\,|x|$ for $|x| \leq \pi/2$ and choosing 
$\lambda = 1/\pi$ in (10.1), this bound yields
\bee
Q(S_n;0) 
 & \leq & 
\Big(\frac{96}{95}\Big)^2 \ \frac{1}{\pi} \int_{-\pi}^{\pi}
\exp\Big\{- 2\sigma^2\,\sin^2(t/2)\Big\}\,dt \\ 
 & \leq & 
\Big(\frac{96}{95}\Big)^2 \ \frac{1}{\pi}  \int_{-\infty}^{\infty}
\exp\Big\{- \frac{2\sigma^2}{\pi^2}\,t^2\Big\}\,dt \ = \ 
\Big(\frac{96}{95}\Big)^2\, \frac{1}{\sigma \sqrt{2}}\,\sqrt{\pi}.
\ene
That is, simplifying the numerical constant, we arrive at
$$
M(S_n) \leq \frac{1.28}{\sqrt{\Var(S_n)}}, \qquad 
\Var(S_n) = \sum_{k=1}^n p_k q_k.
$$

As we now see, the constant $1.28$ in this inequality can be improved 
by virtue of the upper bound (1.4) which yields
\be
M(S_n) \leq \frac{2}{\sqrt{1 + 4\Var(S_n)}} \leq \frac{1}{\sqrt{\Var(S_n)}}.
\en
Moreover, the best universal constant $c>0$ in
\be
M(S_n) \leq \frac{c}{\sqrt{\Var(S_n)}}
\en
is actually better than 1. As was shown in \cite{B-C-V}, the optimal
constant is given by
$$
c \, = \, \max_{\lambda \geq 0}\, \bigg[\sqrt{2\lambda}\,e^{-2\lambda}
\sum_{k=0}^\infty \Big(\frac{\lambda^k}{k!}\Big)^2\bigg] \, \sim \, 0.4688.
$$
One should however mention that the first inequality in (10.2) is
better than (10.3) for small values of $\Var(S_n)$, namely when
$$
\Var(S_n) < \frac{c^2}{4\,(1-c^2)} \, \sim \, 0.0704.
$$
So, it makes sense to consider improvements in the form 
$M(S_n) \leq \Psi(\Var(S_n))$.

Note that the upper bound in (1.4) also provides a similar lower bound
which may equivalently be rewritten as
$$
N_\infty(S_n) \geq \frac 1 4 + \Var(S_n).
$$
In \cite{M-R} this inequality is sharpened and is extended to all 
$\alpha$-entropy powers with $\alpha \geq 2$ as
$$
N_\alpha(S_n) \geq 1 + 2\beta \,\Var(S_n), \quad
\frac{1}{\alpha} + \frac{1}{\beta} = 1.
$$

\vskip5mm
{\bf Acknowledgments.} We would like to thank the referee for
careful reading and valuable comments. The research of S.B. was supported
by the NSF grant DMS-1855575.


\end{document}